\documentclass[twoside,final,a4paper,10pt]{amsart}
\usepackage{times, dsfont}
\title{A few brief comments on the results of manuscript \texttt{arXiv:math/0405153v3}}
\author{Paolo Piccione}
\address{Departamento de Matem\'atica,
Universidade de S\~ao Paulo, Brazil}
\email{piccione@ime.usp.br}
\date{May 22nd, 2007}
\begin{document}
\maketitle
\begin{abstract}
The results in the recently posted manuscript \texttt{arXiv:math/0405153v3} are incorrect.
The correct version of the aimed results is not original. The preprint contains material
from references that are not properly quoted.
\end{abstract}
\bigskip

\noindent \textbf{On the results.}\enspace
In the recently posted paper \texttt{arXiv:math/0405153v3}\footnote{%
This is a revised version (version 3) posted in April 2007 of a preprint posted in May 2004 (version 1), and
that was later withdrawn in October 2004 (version 2).}
it is discussed the computation of the Maslov index
for a symplectic linear system with constant coefficients. The main result, Theorem~1, gives
an incorrect formula for this index, involving the ``signature'' of the upper right block $B$ of the time one flow of the system.
This block does not in general define a symmetric operator, not even under the more restrictive hypotheses assumed
in Theorem~6.3, and thus it does not have a signature. Among the several incorrect arguments, the fatal
mistake seems to arise from a computation on page 15, before the statement of Lemma~6.2, where the author confuses
some identities satisfied by the elements of the \emph{symplectic group} with the identities satisfied by elements
in its \emph{Lie algebra}. The application of Theorem~1 to the case of semi-Riemannian geodesics (Theorem~7.2)
contains another mistake in formula (7.2). Namely, the contribution to the Maslov index given by the initial
endpoint is not equal to the dimension of the manifold, as asserted by the author, but rather by half of
the signature of the metric tensor. This is easily seen using \cite[Remark~4.1, formula~(4.4)]{degenerate}.
Formula (7.2) contains a further minor mistake, originating from an incorrect statement
of Proposition~7.1 (which is elementary, but also incorrectly proven), in that $-JS$ should be considered rather than $JS$.
Strictly speaking, formula (7.2) does not contain one single correct term.
As to Section~5, where the wrong signature formula is not used, the results are incorrect for
several other reasons. For instance, the term containing the sum of $f\left(\frac{\alpha_j}\pi\right)$
in formula~(5.3)
is \emph{not} a Maslov index \emph{relatively to $L_0=\{0\}\oplus\mathds R^n$}, as claimed, but rather to a Lagrangian of the type
$\bigoplus_{l=1}^pL_0^{(l)}$, with summands $L_0^{(l)}$ depending on the decomposition of $H$ into
direct sum $\bigoplus_{l=1}^pH_l$ of normal forms. Thus, also Proposition~5.2 and Corollary~5.3 are false.
\bigskip

\noindent\textbf{On the originality of the results.} \enspace
The core of the manuscript under review is the proof of a formula relating the Conley--Zehnder index of the flow of the system
with the Maslov index (Section~6). The correct version of this formula has already been established  in the very same context and under
the very same assumptions on the final endpoint\footnote{Assumption \textbf{(H1)} in the manuscript is incorrect, too, and it should be replaced by
$\psi(1)L_0\cap L_0=\{0\}$, see \cite[Section~2]{NosPicTau}.}
in \cite[Corollary~2.11]{NosPicTau};
it involves all the four blocks of the time one flow, see \cite[Definition~2.7]{NosPicTau}. In fact, the results in \cite[Section~2]{NosPicTau} are
proved more generally for arbitrary symplectic linear systems, with possibly non constant coefficients.
Relations between the notions of Conley--Zehnder index,
Kashiwara index, Maslov index and the formulas proved in \cite{NosPicTau} have been investigated in several contexts,
and they are well established in the literature, see for instance \cite{gosson}.
The final relation between Conley--Zehnder index, Maslov index and Kashiwara index, formula (6.2),
is a special case of a formula proven in a more general setting in \cite[Proposition~3.27]{degenerate}.
The Maslov index of geodesics in locally symmetric semi-Riemannian manifolds has already been computed
in \cite{JavaPicc}, under no assumption on the endpoints; the formula given in \cite{JavaPicc}
only uses the initial data, and does not assume knowledge of the time one flow.
\bigskip

\noindent\textbf{On the bibliographical references.}\enspace
A correct quotation of reference \cite{abbo} appears in the manuscript; \cite[Chapter~1]{abbo} contains the material
in Section~2 of the paper, as well as the statement and the proof of Proposition~6.1, which provides
the only correct term in the central result of the paper, Theorem~1/Theorem~6.3.
Section~3.1 of the manuscript
is a faithful copy of \cite[Section~3.2, page 19]{degenerate}, although reference \cite{degenerate} does not appear
in the Bibliography of the manuscript. It is easily checked that the original contribution of the author to Section~3.1
concentrates in its final part, starting on line 15 of page 6, which contains entirely incorrect material.
Note that there is no such thing as the differential of $\omega$, appearing in the final centered
formula of Section~3.1; the equality between the first and the last term of such equality is merely
a definition, and it does not prove any of the author's claims. In fact, it is not true that ``the differential
of the map $\varphi_{L_0,L_1}$ \emph{at a point $L$} does not depend on $L_1$'', as asserted.
Section~3.3 of the manuscript, until formula (3.10) on page 9, is a faithful copy of
\cite[Section~3.7]{degenerate}, again, not quoted by the author. Concerning the manuscript
\cite{degenerate} and other related material, the reader should be aware of some authorship claims raised by the author in
a recently posted note (see \cite{Por}).
As to the Maslov index of geodesics, in Section~7, no quotation of reference \cite{JavaPicc} appears
in the manuscript; criticism to the incorrect results in the initial version of the author's manuscript
had already appeared in \cite{JavaPicc}. Instead, on page 17 it is quoted paper \cite{MusPejPor} for a
classical result on the fact that conjugate points do not accumulate at the starting point, while
this result has nothing to do with the authors\footnote{
In fact, the reader will find an incorrect proof of this fact provided by the authors of \cite{MusPejPor}
in the posted version of the manuscript \cite[page 16]{MusPejPor}.}  of \cite{MusPejPor}
(note that in the case under consideration the set of conjugate points is finite, and thus it does not accumulate \emph{anywhere}).
Finally, the author quotes a previous result by
himself as some collaborators of his giving the equality between the Maslov index and ``the total number of conjugate points'',
while it is well known that the Maslov index of a geodesic is \emph{not} equal to the number
of conjugate points. The correct value for the number of conjugate points is computed
in \cite[Corollary~3.5]{JavaPicc}.

\end{document}